\documentclass[a4paper]{article} 

\newif\ifxebiblatex
\xebiblatexfalse

\ifxebiblatex 
\usepackage{fontspec} 
\usepackage{polyglossia,csquotes} 
\setdefaultlanguage{english} 
\setotherlanguage{greek} 
\fi 

\usepackage{amsmath, amssymb, amsthm, amsfonts, bm, mathtools} 
\usepackage[normalem]{ulem}
\usepackage{graphicx, caption, xcolor, placeins} 
\graphicspath{{./figs/}}
\usepackage[skip=6pt]{subcaption}
\usepackage[hidelinks]{hyperref} 
\usepackage{empheq} 
\usepackage{marginnote} 
\theoremstyle{plain}
\theoremstyle{definition}
\theoremstyle{remark}
\newtheorem{remark}{Remark} 

\newcommand{\eps}{\varepsilon}

\usepackage{tikz, pgfplots}
\usetikzlibrary{arrows, scopes, math} 
\pgfplotsset{compat=1.13, 
discont axis/.style={width=1.1\textwidth, height=.75\textwidth, 
enlarge x limits=false, grid=major, legend pos=north east, 
font=\small, 
unbounded coords=jump, 
every axis plot/.style={mark=*, mark size=.1pt}, 
cont plot/.style={mark=none, mark size=2pt}},
cont plot/.style={mark=none, mark size=2pt} 
} 
\usetikzlibrary{arrows.meta, scopes} 
\tikzset{>={Straight Barb[round, angle=60:.12cm 1]}} 

\usetikzlibrary{external} 
\usepackage{stmaryrd} 

\ifxebiblatex 
\usepackage[backend=biber, style=alphabetic, doi=false, isbn=false]{biblatex} 
\fi 

\title{Optimization Methods for One Dimensional Elastodynamics} 
\usepackage{authblk} 

\author[1,2]{Theodoros Katsaounis} 
\author[3]{Grigorios Kounadis} 
\author[3]{Ioanna Mousikou} 
\author[3]{Athanasios E. Tzavaras}
\affil[1]{University of Crete, Heraklion 71409, Greece} 
\affil[2]{Inst.\ of App.\ and Comp.\ Math.\ (IACM), FORTH, Heraklion 71110, Greece} 
\affil[3]{King Abdullah University of Science and Technology, Thuwal 23955-6900, Saudi Arabia} 
\date{~}

\begin{document} 
\maketitle 

\begin{abstract} 
We propose a new approach for solving systems of conservation laws that admit a variational 
formulation of the time-discretized form, and encompasses the p-system or the system of elastodynamics.
The approach consists of using constrained gradient descent for solving an implicit scheme
with variational formulation, while  discontinuous Galerkin finite element methods is used for the spatial
discretization. The resulting optimization scheme performs well, it has an advantage on how it handles oscillations 
near shocks, and a disadvantage in computational cost, which can be partly alleviated by using techniques
on step selection from optimization methods. 
\end{abstract} 




\section{Introduction} 
\label{sec:1} 
The system of elastodynamics is a nonlinear system of hyperbolic conservation laws 
which describes the propagation of longitudinal (or of shear) waves in an elastic medium. 
The same system describes one-dimensional motions  of a gas, and is widely used as 
a paradigm in the theory of conservation laws, then called as the $p$-system. 
It takes the form 
\begin{equation} \label{eq:1} 
\begin{aligned} 
u_t - v_x = 0, \\ 
v_t - \sigma(u)_x = 0, 
\end{aligned} 
\end{equation} 
where $(x,t)\in \mathbb{R}\times \mathbb{R}_+$, and in the elasticity context $u $ is the strain, $v\in \mathbb{R}$ is the velocity, 
and $\sigma(u)$ is a strictly increasing function describing the stress. For longitudinal motions $u > 0$ while for shear motions $u \in \mathbb{R}$.
The system \eqref{eq:1} is supplemented with initial data
\[
u(x,0) = u_0 (x), \quad v(x,0) = v_0 (x).
\]
Smooth solutions of \eqref{eq:1}  generally develop discontinuities in finite time,
reflecting the development of shock waves, and classical solutions cease to exist. 
Introducing the concept of weak solutions, global solutions are constructed using viscosity approximations \cite{Dip83, She94}, 
relaxation approximations \cite{Tza99, Ser00} and numerical schemes (see \cite{GR13} and references therein). 
For the theory of shock waves we refer to \cite{Daf16}. 

Historically, the first methods used to produce numerical approximation of solutions with shocks were finite difference methods, where one
replaces the derivatives of the unknown functions by their finite difference approximations, \cite{Lev92}. Because of the difficulties arising with the application of  finite difference methods to problems with realistic geometries, other methods such as finite volume and finite element methods were later
introduced. Finite volume methods produce approximations for the average of the solution over smaller domains and preserve the conservation property of the exact solution, \cite{Lev02}, \cite{Shu06}. 
In contrast, following a completely different approach, the finite element method produces approximations of the solution as a linear combination of piecewise polynomial functions, \cite{SB11}. Due to stability issues of classical finite element methods for hyperbolic conservation laws, the discontinuous Galerkin method was introduced. The original discontinuous Galerkin method was proposed  by Reed and Hill for solving steady-state neutron transport equations \cite{RH73} and it has been extended for solving nonlinear scalar conservation laws, \cite{CS82}, \cite{CS89C}.
In this framework, the basis functions are completely discontinuous across each element interface and they usually consist of piecewise polynomials defined locally. 
The method was generalized for multidimensional problems, e.g.\ \cite{CS89}, \cite{CS98}, as well as for problems with higher order spatial derivatives which are not necessarily hyperbolic, e.g.\ \cite{CS98A}, \cite{CGL09}.  
Further literature about numerical methods for hyperbolic conservation laws can be found in \cite{GR13}, \cite{Tor13}. 

\medskip 


Our objective is to introduce a new approach for the numerical approximation of \eqref{eq:1} inspired by optimization methods. This approach
is not expected to work for general systems of conservation laws, but applies to special systems that can be viewed as time-discretizations of 
Hamiltonian dynamics and includes in particular \eqref{eq:1}. It is motivated by an approximation developed in 
\cite{DST00} that produces entropy weak solutions. Namely, given a time step $k$, and initial data $(u^0, v^0)$, one produces time iterates $(u^j, v^j)$
by solving the minimization problem
\begin{equation}
\label{eqmin}
\min_{\substack{  \frac{u-u^{j-1} }{k} = v_x }}   \int \tfrac{1}{2} (v - v^{j-1})^2 + W(u) \, dx .
\end{equation}
Here, at the $j$-th  time step, $(u^{j-1}, v^{j-1})$ are given, and $(u^j, v^j)$ is selected as the minimizer of problem \eqref{eqmin}.
For $W(u)$ convex, \eqref{eqmin} consists of minimizing a convex functional over an affine constraint and the iterates are well defined.
Moreover, the iterates yield via interpolation in time approximate solutions that converge as $k \to 0$ to a weak
solution $(u(x,t), v(x,t))$ that decreases all the convex entropies of \eqref{eq:1}, see \cite{DST00}.

In the present work,  the variational scheme \eqref{eqmin} is solved via the method of constrained gradient descent in order
to obtain an explicit scheme. Then we use ideas from finite volumes and from discontinuous Galerkin methods for the spatial discretization, 
utilizing tools from the existing literature, in order to obtain a fully discrete numerical scheme for \eqref{eq:1}. 
We view this as a paradigm to develop and test such optimization
motivated methods for conservation laws with special structure. 
A (nontrivial) variant of the variational scheme \eqref{eqmin} is available  for the system of multi-dimensional elastodynamics with polyconvex energy, see \cite{DST01} 
and \cite{MT13}, and our hope
is to eventually extend this methodology in this interesting context. 
At present, we are interested in developing carefully the one-dimensional case, see Section~\ref{sec:2}. 

The numerical experiments performed here indicate that, when shocks are present in the solution and in the absence of any special techniques that handle oscillations, our method achieves better results compared to some classical DG methods, see Section~\ref{sec:4}. A significant disadvantage however is that the present method is computationally demanding; to mitigate this issue, we propose some novel computational techniques and mechanisms that significantly improve the convergence speed. 
We are working in implementing the scheme in multiple dimensions, where we expect further benefits. 

The manuscript is organized as follows. In Section~\ref{sec:2} we briefly describe the derivation of two variational schemes that arise from time discretization of the original system. 
In Section~\ref{sec:3} we introduce the constrained gradient descent, the optimization method that will be used to solve the minimization problem.
Then we formally describe the discontinuous finite element space, and state the fully discrete form of the method. We also give the required implementation details concerning various aspects of the method. 
Finally, in Section~\ref{sec:4} we verify numerically that the rate of convergence is optimal for smooth solutions; we investigate the quality of solutions in the presence of shocks, and compare the total variation with the classical RKDG discretization of the system as conservation law. 
In the last part of the section we assess the speed of gradient descent and propose techniques to accelerate the convergence. 


\section{The approximation framework} 
\label{sec:2} 

We work in a bounded domain, $(x,t) \in [0,1]\times[0,T]$ with periodic boundary conditions. 
Let  $k$ be a time step, $\{t^j\}_{j=0}^M$ be a partition of $[0,T]$ with  time step $k$, and we are interested in constructing an approximate 
solution of \eqref{eq:1} as follows. We start with the initial data $(u^0, v^0)$ periodic and of zero mean. At each time step $t^j= j k$ given $(u^{j-1}, v^{j-1})$ we solve the
implicit problem
\begin{equation} \label{eq:2} 
\begin{aligned} 
\frac{u^j-u^{j-1}}{k} &= v_x^j, \\ 
\frac{v^j-v^{j-1}}{k} &= \sigma(u^j)_x.     
\end{aligned} 
\end{equation} 
The solutions $(u^j (x), v^j (x))$ of problem \eqref{eq:2} are obtained by variational minimization.

We give below a brief outline of two variational schemes developed in \cite{DST00} to solve \eqref{eq:2}  
and describe the limit  to \eqref{eq:1}  as $k \to 0$. The reader is referred to \cite{DST00} for details. 
For the first scheme, we work with the equivalent second order equation 
\begin{equation} \label{eq:3} 
y_{tt} - (\sigma(y_x))_x = 0, 
\end{equation} 
which, by replacing $u=y_x$ and $v=y_t$, reduces to \eqref{eq:1}. We will assume that $y$ is periodic in space, and discretize \eqref{eq:3} in time by 
\begin{equation} \label{eq:4} 
\frac{y^{j} - 2y^{j-1} + y^{j-2}}{k^2} - \bigl(\sigma(y_x^{j})\bigr)_x = 0. 
\end{equation} 
The variational scheme consists of, given $y^{j-1}$, $y^{j-2}$ functions of zero-mean, to find the minimizer  $y^{j}$
of the problem 
\begin{equation} \label{eq:5} 
\min_{y\in H^1_\text{per}} I_{k,j}[y] = \min_{y\in H^1_\text{per}}\int W(y_x) + \frac{ \big ( y-2y^{j-1} + y^{j-2} \big )^2 }{2k^2}\,\mathrm dx,  
\end{equation} 
where the minimization is performed over $H^1_\text{per}$, the periodic $H^1$ functions with zero mean, and $W'(u) = \sigma(u)$. Under some growth and smoothness assumptions for $W(u)$,  which is principally assumed as strictly convex,
it is shown that $I_{k,j}$ attains a unique minimizer, say $y^{j}$, 
and the minimizer satisfies the Euler-Lagrange  equations 
\begin{equation} \label{eq:6} 
\int \sigma(y_x^{j})\zeta_x + \frac{ y^{j} - 2y^{j-1} + y^{j-2} }{k^2}\zeta\,\mathrm dx,\quad \forall \zeta\in W^{1,p}, 
\end{equation} 
stating that \eqref{eq:4} holds in a weak sense. 

Letting now 
\begin{equation}
\label{transf}
u^{j} = y_x^{j},\quad v^{j} = \frac{y^{j}-y^{j-1}}{k},
\end{equation}
a second variational scheme is derived as follows. First, note that the transformation \eqref{transf} reduces
equation \eqref{eq:4} into the equivalent system \eqref{eq:2}. Second, as shown below,   \eqref{eq:2}
are the Euler-Lagrange equations of the constrained minimization problem:
Given $(u^{j-1},v^{j-1})$, find $(u^{j}, v^{j})$ 
the minimizer of  the problem
\begin{subequations} \label{eq:7} 
\begin{equation} \label{eq:7a} 
\min J [u,v] = \min\int W(u) + \frac{(v-v^{j-1})^2}{2}\,\mathrm dx, 
\end{equation} 
where the minimization is done over the set of functions satisfying the affine constraint
\begin{equation} \label{eq:7b} 
\int \frac{u-u^{j-1}}{k}\phi + v\phi_x\,\mathrm dx = 0,\quad \forall \phi\in C^1. 
\end{equation} 
\end{subequations} 
The problem  \eqref{eq:7a}-\eqref{eq:7b}  attains a unique minimum at, say, $(u,v)$, see \cite{DST00}.

The Euler-Lagrange equations for \eqref{eq:7a}-\eqref{eq:7b} are computed as follows:
 Let $(u,v)$ be the minimizer and consider a variation
$(u + \delta p, v + \delta q)$ where $p, q$ are smooth test functions.
Then \eqref{eq:7a} implies
\[
\int W(u + \delta p) + \frac{(v + \delta q - v^{j-1})^2}{2}\,\mathrm dx \ge \int W(u) + \frac{(v-v^{j-1})^2}{2}\,\mathrm dx, \quad  \forall \delta \in \mathbb{R}. 
\]
Taking the limits $\delta  > 0, \; \delta \to 0$ and then $\delta < 0, \; \delta \to 0$ we deduce
\begin{equation} \label{eq:8} 
\int\sigma(u)p+(v-v^{j-1})q\,\mathrm dx = 0,\quad \forall p, q,\ \text{smooth}.
\end{equation} 
The constraint \eqref{eq:7b} implies
\[
\int \frac{u + \delta p -u^{j-1}}{k}\phi + (v + \delta q) \phi_x\,\mathrm dx = 0,\quad \forall \phi\in C^1, 
\]
which, taking $\delta \to 0$ implies that the variations satisfy $p = k q_x$. Hence, the minimizer satisfies the Euler-Lagrange 
equations \eqref{eq:7b} and 
\begin{equation} \label{eq:9} 
\int  (v-v^{j-1})q   +  k \sigma(u) q_x \,\mathrm dx = 0,\quad \text{for}\ q\ \text{smooth}, 
\end{equation} 

The minimizing scheme produces iterates $(u^j (x), v^j (x))$, with $j = 1, ... , M k$, $M k =T$. 
Using the iterates we define
approximate solutions $(u^k (x, t), v^k (x,t))$ to \eqref{eq:1} via either piecewise constant, or piecewise linear
interpolation in time. 
It is shown in \cite{DST00} using the theory of compensated compactness
that $(u^k, v^k)$ converge to $(u,v)$ almost everywhere (both piecewise constant and piecewise linear interpolations yield the same limit)
and that $(u,v)$ satisfies \eqref{eq:1} and the following form of entropy inequalities:
For any entropy pair $(\eta, q)$ such that $\nabla\eta \nabla f = \nabla q$ with $f(u,v) = (-v, -\sigma(u))^\mathsf T$ the solution $(u(x,t), v(x,t))$
satisfies
\begin{equation}
\label{entradm}
\partial_t \eta(u,v) + \partial_x q(u,v)\leq 0, 
\end{equation}
in distributions for any  entropy $\eta(u,v)$ convex. Observe that \eqref{entradm} is the same admissibility condition that artificial viscosity
would produce for the system \eqref{eq:1}.


\section{Numerical method} 
\label{sec:3} 
In this section we describe the numerical scheme that we will utilize. We solve the minimization problem by a constrained gradient descent method.
First we describe the method in an abstract framework.
Subsequently, the method is adapted to the problem at hand, \eqref{eq:7a}-\eqref{eq:7b},  and the
resulting  weak formulation is expressed in a continuous Finite Element (FE) space. 
To approximate solutions containing shocks we then move to a discontinuous FE space and introduce some necessary stabilization terms. 
In \ref{subsec:3p3} we state the fully discrete formulation of the method, assess the computational complexity and various other aspects. 
Finally, in the last subsection we list important details concerning the implementation.

\subsection{Constrained Gradient Descent}
\label{subsec:3p}

Consider the constrained minimization problem 
\begin{subequations} 
\begin{equation} \label{eq:min} 
\min_{x \in \mathcal{A}}  F(x)
\end{equation} 
where $X$ is a Banach space, $F : X \to \mathbb{R}$ is a convex functional, and the minimization is done over an affine subspace $\mathcal{A}$
\begin{equation} \label{eq:constr} 
\mathcal{A} = \big \{ x \in X :   A [x] = c \big \} 
\end{equation} 
\end{subequations} 
defined by linear functionals $A = (A_1 , .... , A_n)  : X \to \mathbb{R}^n$ with $c \in \mathbb{R}^n$. 
This problem consists of minimizing a convex function $F$ over an affine subspace. Under fairly general conditions: $X$ is reflexive, $F$ is convex,
coercive and weakly lower semicontinuous on $X$, while the linear functionals $A_i$ determine a weakly closed subspace of $X$, this 
minimization problem has a solution \cite[Cor. 3.23]{Bre83}.
Moreover, when $F$ is strictly convex the solution is unique. The associated Euler-Lagrange equations define the
minimizer $x$ implicitly. 

We propose to compute the minimizer via gradient descent taking also into account the affine constraint \eqref{eq:constr}. 
(The method is expected to work when the constraint is affine, and it would lead in general to nonconvex problems when
the constraint is nonlinear.) Given an iteration step $\lambda$ and $x_l \in \mathcal{A}$, the gradient descent method computes the 
next iterate $x_{l+1}$ by 
\begin{equation}
\label{eq:cgd}
\begin{aligned}
x_{l+1} - x_l &= \lambda \frac{\delta F}{\delta x} (x_l) 
\\ 
x_{l+1} &\in \mathcal{A} \, .
\end{aligned}
\end{equation}
The variational derivative $ \frac{\delta F}{\delta x}$ for the constrained problem \eqref{eq:min} is computed by
\begin{equation}
\label{eq:vd}
\big \langle  \frac{\delta F}{\delta x}(x) , \varphi \big \rangle := \lim_{\substack{ \eps \to 0 \\[2pt]  x \in \mathcal{A}   \\[2pt]  x + \eps \varphi \in \mathcal{A}   }}   \frac{ F( x + \eps \varphi) - F(x) }{\eps} \, .
\end{equation}
where $\varphi$ is a test function. Equation \eqref{eq:vd} precisely defines the variational derivative of the constrained minimization problem \eqref{eq:min}-\eqref{eq:constr}. In applications it
will be expressed by introducing a basis function on the constraint subspace $\mathcal{A}$.

\subsection{Adaptation to the specific minimization problem } 
\label{subsec:3p1} 
Given $(u^{j-1} , v^{j-1})$, the  $j$-th iterate $(u^j, v^j)$ is constructed as the solution of the constrained minimization problem 
\begin{subequations} 
\begin{equation} \label{eq:10a} 
\min_{u,v} J[u,v] = \min_{u,v} \int\biggl(\frac{(v-v^{j-1})^2}{2} + W(u)\biggr)\,\mathrm dx 
\end{equation} 
subject to the affine constraint
\begin{equation} \label{eq:10b} 
\frac{u-u^{j-1}}{k} = v_x, 
\end{equation} 
\end{subequations} 
where $W'(u) = \sigma(u)$. 

To implement gradient descend (GD) we first have to calculate the variational derivative of \eqref{eq:10a}. If we let $\varepsilon P := (\varepsilon p, \varepsilon q)^\mathsf T$ be the variation of $U := (u, v)^\mathsf T$ in the direction of $(p,q)^\mathsf T$, and $j(\varepsilon) := J(U + \varepsilon P)$, the derivative is
\begin{equation} \label{eq:11} 
j'_P(0) =  \Big \langle \frac{\delta J}{\delta U}, (p, q)  \Big \rangle = \int (v-v^{j-1})\, q  + W'(u)\, p\,\mathrm dx. 
\end{equation} 
From the constraint \eqref{eq:10b} we have that $p = k\,q_x$, 
therefore by substituting in \eqref{eq:11} we get 
\[
j'_P(0) = \int (v-v^{j-1})\,q + W'(u)\,k\,q_x\,\mathrm dx.   
\]
Recalling that we use as template periodic boundary conditions, 
we work in the Sobolev space $H^1_\text{per}(0,1)$. Given some approximation $u_l$, $v_l$ to $u$, $v$, the gradient descent method (GD) will decrease the value of \eqref{eq:10a} by finding $v_{l+1}\in H^1_\text{per}(0,1)$ such that 
\begin{subequations} \label{eq:12} 
\begin{equation} \label{eq:12a} 
(v_{l+1}, \phi) = (v_l,\phi) - \lambda \bigl( (v_l-v^{j-1},\phi) + k(W'(u_l),\phi_x) \bigr),\quad \forall \phi\in H^1_\text{per}(0,1), 
\end{equation} 
where $\lambda$ is the GD iteration step.
Let  $A: H^1_\text{per}\times H^1_\text{per}\to \mathbb R$ stand for the bilinear form $A(v,\phi)=(v,\phi)$, and let $G_1 (\phi) = G_1 (\phi; \lambda, u_l, v_l, v^{j-1})$  
denote the right hand side  of \eqref{eq:12a}, which is then expressed as  $A(v_{l+1},\phi)= G_1 (\phi)$. 

The constraint \eqref{eq:10b} is enforced by  defining $u_{l+1}\in H^1_\text{per}(0,1)$ via
\begin{equation} \label{eq:12b} 
(u_{l+1},\phi) = (u^{j-1},\phi) + k(v_{x,l+1},\phi),\quad \forall \phi\in H^1_\text{per}(0,1) \, .
\end{equation} 
\end{subequations} 
If we set $G_2 (\phi; v_{l+1}, u^{j-1})$ to be the right hand side then \eqref{eq:12b} is expressed via
$A(u_{l+1},\phi) = G_2 \phi; v_{l+1}, u^{j-1})$. 

We have implemented the aforementioned scheme, \eqref{eq:12}, and verified numerically that the order of convergence is optimal. However, continuous finite element spaces are not suitable for solutions containing shocks, since spurious oscillations of large amplitude are formed near discontinuities. 
In fact it can be shown that the standard FE discretization of the original system as a conservation law using piecewise linear functions reduces to a central difference method that is unstable, see \cite[Ch.~10]{Lev3}.

\subsection{Discretization} 
\label{subsec:3p2} 
We will work in the well established framework of Discontinuous Galerkin methods (DG). The original system is a conservation law; conservation laws, in the setting of DG, have been studied by Cockburn and Shu in the excellent series of papers \cite{CS89, CS89A}, from which we will borrow various tools. 

Let $\{x_{i+\frac 1 2}\}_{i=0}^N$ be a partition of $[0,1]$,  let 
$I_i$ denote the cell $\bigl[x_{i-\frac 1 2}, x_{i+\frac 1 2}\bigr]$ and  let $h_i$ be the its length $h_i=x_{i+\frac 1 2} - x_{i-\frac 1 2}$. We seek a solution in the space of piecewise polynomial functions of order $K$ with periodic boundary conditions 
\[ 
V_h = \Bigl\{\phi\in L^1(0,1):\ \phi\big|_{I_i}\in\mathbb P^K(I_i),\ i=1,\ldots,N,\ \phi(0)=\phi(1)\Bigr\}.
\] 
Due to their various beneficial properties, and as commonplace in DG methods, we use Legendre polynomials as basis functions. Let $\{\phi_i^\ell\}_{\substack{i=1,\ldots, N\\ \ell=0,\ldots,K}}$ be the basis, where $\phi_i^\ell$ is the $\ell$-th degree Legendre polynomial scaled onto the cell $I_i$. 

Having defined the general GD step in \eqref{eq:12} we now proceed in finding $v^j$, $u^j$, the solution at the next time step, in the setting of the now discontinuous finite element space $V_h$. 
Let $\{u_l^j\}_{l=0}^{L_j}$, $\{v_l^j\}_{l=0}^{L_j}$ represent the sequence generated by GD to approximate $u^j$, $v^j$. 
A sensible initial guess, to kickstart the method, is the solution in the previous time step, $v^j_0 = v^{j-1}$, $u^j_0 = u^{j-1}$; (in the first time step we will use the $L^2$-projection to approximate the initial data, $u^0(x) = \operatorname{P} u(x,0)$, $v^0(x) = \operatorname{P} v(x,0)$).  
Having calculated the first $l$ iterates, $v_{l+1}^j$ is the unique function in $V_h$ that satisfies 

\begin{subequations} \label{eq:13} 
\begin{equation} \label{eq:13a} 
A(v_{l+1}^j,\phi) + \int_{\Gamma} \frac\mu h\, \llbracket v_{l+1}\rrbracket\cdot\llbracket\phi\rrbracket\,\mathrm ds = G_1 (\phi; u_l^j, v_l^j,v^{j-1}),\quad \forall\phi\in V_h, 
\end{equation} 
where $\Gamma$ is the boundary of all elements (the cell interfaces in one dimension), $h$ the minimum cell length, $\mu$ a positive constant (the penalty), and $\llbracket u\rrbracket:= u^+ - u^-$ is the classical jump operator where $u^+, \ u^-$ denote solution to the right and left of the interface, respectively. 
A term that penalizes jumps across cell interfaces is required; this is motivated by the theory of DG methods for elliptic equations, see \cite{ABCM02}, where it plays the role of stabilization term; for given $\mu$ large enough, it ensures that the corresponding bilinear form is coercive in a suitable norm. 
We enforce the constraint by updating $u_{l+1}^j$ as 
\begin{equation} \label{eq:13b} 
A(u_{l+1}^j, \phi) = G_2 (\phi; v_{l+1}^j,u^{j-1}),\quad \forall\phi\in V_h. 
\end{equation} 
\end{subequations} 
The above procedure is repeated until the solution converges and the integral  to be minimized, 
\begin{equation} \label{eq:14} 
I_l = I[v_l,u_l] := \int_0^1 \left(\frac{(v_l^j-v^{j-1})^2}{2} + W(u_l)\right)\mathrm dx, 
\end{equation} 
stops decreasing; more details about the stop criteria will follow.

\subsection{Solution of the discretized problem} 
\label{subsec:3p3} 
Let $u_h, v_h\in V_h$ be the numerical solution at $j$-th time step. Taking into consideration that the support of $\phi_i^\ell$ is $I_i$, from \eqref{eq:13a}, $v_h$ is the unique function that satisfies 
\begin{subequations} \label{eq:15} 
\begin{multline} \label{eq:15a}
\int_{I_i}v_{h,l+1}\phi_i^\ell\,\mathrm dx + \left.\frac \mu h \llbracket v_{h,l+1}\rrbracket\llbracket\phi_i^\ell\rrbracket\right|_{x_{i-\frac 1 2}}^{x_{i + \frac 1 2}} = \int_{I_i} v_{h,l}\phi_i^\ell\,\mathrm dx  -\\ 
 \lambda\int_{I_i} (v_{h,l}-v_h^{j-1})\phi_i^\ell - k\,W'(u_{h,l})\phi_{x,i}^\ell\,\mathrm dx, \quad \forall i,\ \ell, 
\end{multline} 
similarly from \eqref{eq:13b}, by integrating by parts, $u_h$ is given by 
\begin{equation} \label{eq:15b} 
\int_{I_i} u_{h,l+1}\phi_i^\ell\,\mathrm dx = \int_{I_i} u_h^{j-1}\phi_i^\ell\,\mathrm dx - k\int_{I_i} v_{h,l+1}\phi_{x,i}^\ell\,\mathrm dx + k\, \widehat v_{h,\ell+1}\phi_i^\ell\Bigr|_{x_{i-\frac 1 2}^+}^{x_{i + \frac 1 2}^-}, \quad \forall i,\ \ell, 
\end{equation} 
\end{subequations} where the last term arises from integration by parts; $\widehat v_h$ is the {\em numerical flux}, it comes from the FV literature where the system flux is approximated by solving a Riemann problem on cell interfaces, and is a classical way to stabilize the DG form. More information about the numerical flux selection is given in the next section. 

\noindent
If we define the matrices 
\begin{align*} 
&M_{m,n} = M_{iK+\ell, jK+\ell'} = (\phi_i^\ell, \phi_j^{\ell'}), & 
&S_m = S_{jK+\ell} = (W'(u_{h,l}), \phi_{x,j}^\ell), \\ 
&A_{m,n} = A_{iK+\ell, jK+\ell'} = (\phi_i^\ell, \phi_{x,j}^{\ell'}), & 
&J_{m,n} = J_{iK+\ell, jK+\ell'} = \llbracket \phi_i^\ell\rrbracket\llbracket\phi_j^{\ell'}\rrbracket, 
\end{align*} 
equations \eqref{eq:15} may be written as 
\begin{align*} 
(M+J)\, D_{l+1} &= M\bigl(D_l - \lambda (D_l-D^{j-1})\bigr) - \lambda kS_l, \\ 
M\,C_{l+1} &= M C^{j-1} - kA D_{l+1} + kB_l, 
\end{align*} 
where the vectors $C$ and $D$ hold the degrees of freedom for $u_h$ and $v_h$ respectively and $B$ holds the values arising from the numerical flux. We observe that all matrices are block diagonal; furthermore, due to the orthogonality of Legendre polynomials, the mass matrix $M$ reduces to a diagonal $M_{iK+\ell, jK+\ell'} = \frac{2}{2\ell+1}\delta_{\ell\ell'}$ (for $i=j$) and the matrix $A$ reduces to a strictly lower triangular where the nonzero elements are $A_{iK+\ell, jK+\ell'} = 2$ for $\ell=\ell'{-}1, \ell'{-}3, \ldots$ (and $i=j$). Finally, numerical flux calculations are simplified by the fact that $\phi_i^\ell(x_{i+\frac 1 2}^-)=1$ and $\phi_i^\ell(x_{i+\frac 1 2}^+)=(-1)^\ell$. 

\begin{remark} 
It is easy to see that \eqref{eq:15} is consistent with the continuous problem \eqref{eq:12}. In particular, if solution $(u,v)$ is continuous, the penalty term in \eqref{eq:15a} vanishes while the numerical flux in \eqref{eq:15b} reduces to the original flux. After integration by parts \eqref{eq:15b} coincides with \eqref{eq:12b}. 
\end{remark} 

Most of these are well known properties of DG methods that reduce the computational complexity and improve the parallelizability of the code. 
Compared to continuous FE methods, no quasiuniformity of the grid is required. 
In addition, $h$- and $p$-adaptivity and the handling of complex geometries (in the case of multiple dimensions) are made significantly simpler.

\subsection{Implementation details} 
\label{subsec:3p4} 

\paragraph{Numerical flux.} To handle the term that arises by the integration by parts in \eqref{eq:15b} we will employ the Local Lax Friedrich (LLF) numerical flux. For a general system of conservation laws, $U_t + f(U)_x = 0$, the LLF numerical flux is  
\[ 
\widehat{f}^\text{LLF}(U^-,U^+) = \tfrac 1 2 \bigl(f(U^-) + f(U^+) - \alpha(U^+-U^-)\bigr), 
\] 
where $\alpha=\max_{(\min(U^-,U^+), \max(U^-,U^+))} |\lambda^{\max}(U)|$ and $\lambda^{\max}$ the maximum eigenvalue of the Jacobian $\frac{\partial f}{\partial U}$.
This flux is simple to implement and computationally efficient, but introduces a significant amount of numerical diffusion. It is worth noting that the choice of the numerical flux does not have a significant impact as the polynomial degree increases. We refer to \cite{Tor13} for a study about numerical fluxes. 

\noindent
In our system the numerical flux takes the form  
\[ 
\widehat v(x_{i+\frac 1 2}) = \frac 1 2 \Bigl( v_h(x_{i+\frac 1 2}^+)+v_h(x_{i+\frac 1 2}^-) - \alpha\bigl(u_h(x_{i+\frac 1 2}^+)-u_h(x_{i+\frac 1 2}^-)\bigr)\Bigr),  
\] 
with $\lambda^{\max} = \sqrt{\sigma'(u)}$. 

\paragraph{Slope limiter.} The presence of spurious oscillations near discontinuities is a known phenomenon when applying  DG methods for conservation laws; and the jump penalty term we had to introduce in \eqref{eq:15a} does not improve the situation. A technique to combat this phenomenon is limiting the slope of the solution based on its value in adjacent cells. 
The most well known limiter, colloquially known as {\itshape minmod}, is presented at \cite{CS89}. In short; in each cell we write the solution at the endpoints as 
\[ 
u_{i+\frac 1 2}^- = u_i^0 + \tilde u_i,\quad u_{i-1/2}^+ = u_i^0-\tilde{\tilde{u}}_i, 
\] 
we modify the solution by limiting $\tilde u_i$, $\tilde{\tilde{u}}_i$ by 
\begin{equation} \label{eq:16} 
\tilde u_i^\text{mod} = m(\tilde u_i,\, u_{i+1}^0{-}u_i^0,\, u_i^0{}-u_{i-1}^0),\ \ \tilde{\tilde{u}}_i^\text{mod} = m(\tilde{\tilde{u}}_i,\, u_{i+1}^0{-}u_i^0,\, u_i^0{-}u_{i-1}^0), 
\end{equation} 
where $m$ is the modified minmod function.
The degrees of freedom can be then calculated using \eqref{eq:16}. To ensure that the solution is TVD, this procedure has to be performed after projecting the solution to the characteristic fields of the system.

The degrees of freedom can be uniquely determined from \eqref{eq:16} for up to quadratic polynomials. For higher order polynomials, the usual procedure is to set $u_i^\ell=0$ for $\ell=2, \ldots, K$ in the cells where the limiter is applied, which effectively limits the approximation order and thus the accuracy. To overcome this issue {\itshape moments limiter} is introduced in \cite{BDF94}. We iteratively limit the degrees of freedom by 
\begin{equation} \label{eq:17}  
(2\ell+1)u_i^{\ell+1} = m((2\ell+1)u_i^{\ell+1},\, u_{i+1}^\ell{-}u_i^\ell,\, u_i^\ell{-}u_{i-1}^\ell), 
\end{equation} 
starting from the highest degree, $\ell=K-1$, and moving down until we encounter an $\ell$ for which the degree of freedom is not modified by \eqref{eq:17}. Moments limiter successfully maintains the order of accuracy most of the times. 

We will use minmod limiter for up to quadratic polynomials and moments limiter for higher order ones. We limit both $v_h$ and $u_h$ after each time step. Applying the slope limiter after each GD step significantly increases the computational complexity without having appreciable benefits. 

\paragraph{Gradient Descent.} The number of iterations of GD play a significant role in the computational complexity of the method. We control the convergence of the algorithm using three quantities: 
\begin{itemize} \setlength\itemsep{0pt} 
\item[-] the difference between two successive evaluations of the integral-to-be-minimized \eqref{eq:14}, $|I_{l+1}-I_l| < c_I$, with default tolerance $c_I = 10^{-14}$; 
\item[-]  the difference of two successive approximations of $u$, $\|u_{l+1}-u_l\|<c_u$, with default tolerance $c_u=10^{-14}$; 
\item[-]  the number of actual iterations performed, $c_i$, with default value $c_i=250$. 
\end{itemize} 
The method is said to have converged when both conditions $|I_{l+1}-L_l|<c_I$ and $\|u_{l+1}-u_l\|<c_u$, are satisfied or when the maximum number of iterations is reached, where $\|\cdot\|$ denotes here and throughout the $L^2$ norm. 

The choice of the GD step, $\lambda$, is also crucial. A very large value might prevent the method from converging, while a too small value might require a large amount of iterations. The value of step $\lambda$  can also be controlled adaptively.  A simple heuristic algorithm that improves significantly the convergence speed is presented in Section~\ref{subsec:4p3}; alternative algorithms  for choosing $\lambda$, can be found in the optimization literature, see for example \cite{KW19}. 


\section{Numerical results} 
\label{sec:4} 
In this section we present the results of a series of numerical experiments demonstrating the approximating features, effectiveness and robustness of the method. We will use $\sigma(u)=u^3+u$, and set the penalty constant $\mu=1$. In the first parts of this section we are interested in the approximation properties of the method, therefore we will use a fixed $\lambda = 1/4$ and some --rather strict-- GD stop criteria (described in the previous section).

\subsection{Effective order of convergence} 
\label{subsec:4p1} 
To verify the convergence rate of the method we consider the following smooth initial conditions  
\begin{equation} \label{eq:18} 
u(x,0) = u_0(x) = 2-\exp(-0.5(x-4)^4),\quad v(x,0) = v_0(x) = u_0'(x), 
\end{equation} 
for $x\in[0,8]$ and $T=1/40$. 
The solution remains smooth for the duration of the simulation; a separate algorithm (that discretizes the original system as a conservation law) has been used to obtain a high quality approximation of the solution for error estimation purposes. 
We consider piecewise linear polynomials and a uniform partition in space. The spatial rate of convergence of the method is expected to be $2$, while the temporal only $1$. To factor in  the spatial accuracy in our calculations we set $k=c_\text{cfl}/h^2$, where $c_\text{cfl} = c_\text{RK} \max_x\sqrt{\sigma'(u)}$; the maximum eigenvalue is evaluated at each time step and $c_\text{RK} = 1/8$ is a constant that depends on the degree of polynomials used. Convergence rates are as expected, and can be seen in Table~\ref{tab:1}. 

\begin{table}[htbp] \ttfamily \small \centering 
\setlength{\tabcolsep}{4pt}
\begin{tabular}{|c||c|c|c|c||c|c|c|c|} \hline 
$N$ & $\|u-u_\text{ex}\|$ & {\rmfamily rate} & $\|u-u_\text{ex}\|_\infty$ & {\rmfamily rate} & $\|v-v_\text{ex}\|$ & {\rmfamily rate} & $\|v-v_\text{ex}\|_\infty$ & {\rmfamily rate} \\ \hline 
 20 & 3.166e-02 & -    & 4.146e-02 & -    & 1.061e-01 & -    & 1.134e-01 & -    \\ \hline 
 40 & 1.012e-02 & 1.66 & 1.639e-02 & 1.34 & 2.336e-02 & 2.18 & 3.068e-02 & 1.89 \\ \hline 
 80 & 2.565e-03 & 1.98 & 4.245e-03 & 1.95 & 5.330e-03 & 2.13 & 8.012e-03 & 1.94 \\ \hline 
160 & 6.279e-04 & 2.03 & 1.043e-03 & 2.02 & 1.299e-03 & 2.07 & 2.049e-03 & 1.97 \\ \hline 
320 & 1.552e-04 & 2.02 & 2.595e-04 & 2.01 & 3.225e-04 & 2.01 & 5.149e-04 & 1.99 \\ \hline 
\end{tabular} 
\caption{Discontinuous Galerkin-Interior Penalty (IP), LLF flux, convergence rates. $T=0.025$. Smooth solution. Linear polynomials. \label{tab:1}}
\end{table}

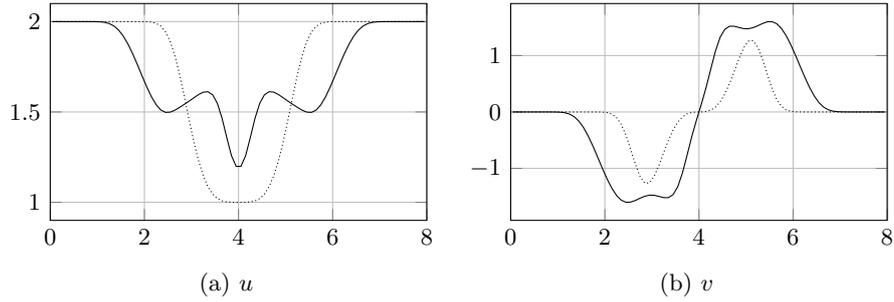
\begin{figure}[htbp] 
\centering
\begin{subfigure}{.49\textwidth} 
\begin{tikzpicture} 
\begin{axis}[discont axis, xmin=0, xmax=8]
\addplot[cont plot, densely dotted] table[x=xn, y=u0] {figs/smooth_linear.txt};\label{fig:1_ic}
\addplot[cont plot] table[x=xn, y=uJ] {figs/smooth_linear.txt};\label{fig:1_T}
\end{axis}
\end{tikzpicture}
\caption{$u$}
\end{subfigure} 
\begin{subfigure}{.49\textwidth} 
\begin{tikzpicture} 
\begin{axis}[discont axis, xmin=0, xmax=8]
\addplot[cont plot, densely dotted] table[x=xn, y=v0] {figs/smooth_linear.txt}; 
\addplot[cont plot] table[x=xn, y=vJ] {figs/smooth_linear.txt};
\end{axis}
\end{tikzpicture}
\caption{$v$} 
\end{subfigure}
\caption{DG-IP, $N=80$, p.w.\ polynomials. Smooth solution. Dotted line (\ref{fig:1_ic}): initial condition, solid line (\ref{fig:1_T}): solution at time $T=0.25$. \label{fig:1} 
} 
\end{figure}

\noindent
To verify the validity of the numerical solution we plot its evolution at time $T=1/4$, along with the initial conditions, in Figure~\ref{fig:1}. 

\subsection{Evolution of discontinuous initial profile} 
\label{subsec:4p2} 
In case of discontinuous solution, and in the absence of any special treatment, i.e.\ slope limiting, some spurious oscillations are generated near discontinuities. An analogous phenomenon is also observed in the standard discretization of the system as a conservation law using DG methods. Oscillations grow larger as the rate $k/h$ gets smaller. 

\begin{figure}[htbp] 
\centering
\begin{subfigure}{.49\textwidth} 
\begin{tikzpicture} 
\begin{axis}[discont axis]
\addplot[cont plot, gray, densely dotted] table[x=xn, y=u] {figs/discont_exact.txt};\label{fig:2_ex}
\addplot[] table[x=x, y=u1] {figs/discont_linear.txt};\label{fig:2_num} 
\end{axis}
\end{tikzpicture}
\caption{$u$}
\end{subfigure} 
\begin{subfigure}{.49\textwidth} 
\begin{tikzpicture} 
\begin{axis}[discont axis]
\addplot[cont plot, gray, densely dotted] table[x=xn, y=v] {figs/discont_exact.txt};
\addplot[] table[x=x, y=v1] {figs/discont_linear.txt}; 
\end{axis}
\end{tikzpicture}
\caption{$v$} 
\end{subfigure}
\caption{DG-IP, $N=160$, $T=0.25$, $k/h=1/8$, p.w.\ linear polynomials. Discontinuous sol. Solid line (\ref{fig:2_num}): num.\ sol., dotted line (\ref{fig:2_ex}): exact sol.\label{fig:2}} 
\vspace{-2ex} 
\end{figure}
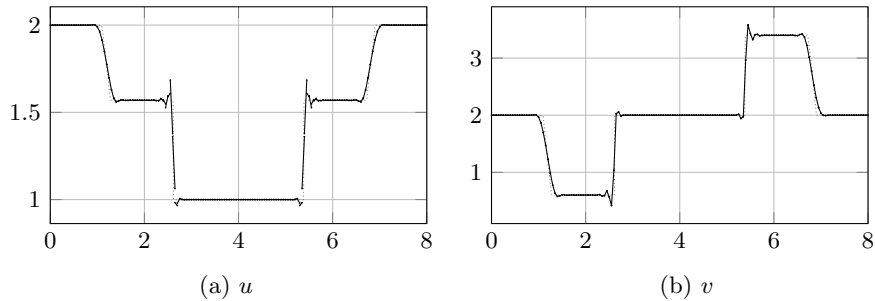

Consider the following discontinuous initial conditions 
\begin{equation} \label{eq:19} 
u_0(x) = \begin{cases} 
1,& \text{if}\ 4\leq x\leq 6, \\ 
2,& \text{otherwise}, 
\end{cases} \qquad v_0(x) = 2. 
\end{equation} 
The solution for piecewise linear polynomials and $k/h=1/8$ at time $T=1/4$ can be seen in Figure~\ref{fig:2}, where the presence of oscillations is evident. 

\begin{figure}[htbp] 
\centering
\begin{subfigure}{.49\textwidth} 
\begin{tikzpicture} 
\begin{axis}[discont axis]
\addplot[cont plot, gray, densely dotted] table[x=xn, y=u] {figs/discont_exact.txt};
\addplot[mark repeat=11, mark phase=10] table[x=x, y=u3] {figs/discont_cubic.txt}; 
\addplot[only marks, mark repeat=11] table[x=x, y=u3] {figs/discont_cubic.txt};
\end{axis}
\end{tikzpicture}
\caption{$u$}
\end{subfigure} 
\begin{subfigure}{.49\textwidth} 
\begin{tikzpicture} 
\begin{axis}[discont axis]
\addplot[cont plot, gray, densely dotted] table[x=xn, y=v] {figs/discont_exact.txt};
\addplot[mark repeat=11, mark phase=10] table[x=x, y=v3] {figs/discont_cubic.txt}; 
\addplot[only marks, mark repeat=11] table[x=x, y=v3] {figs/discont_cubic.txt};
\end{axis}
\end{tikzpicture}
\caption{$v$} 
\end{subfigure}
\caption{DG-IP, $N=160$, $T=0.25$, $k/h=1/28$, piecewise cubic polynomials. Discontinuous solution. \label{fig:3}} 
\vspace{-2ex} 
\end{figure}
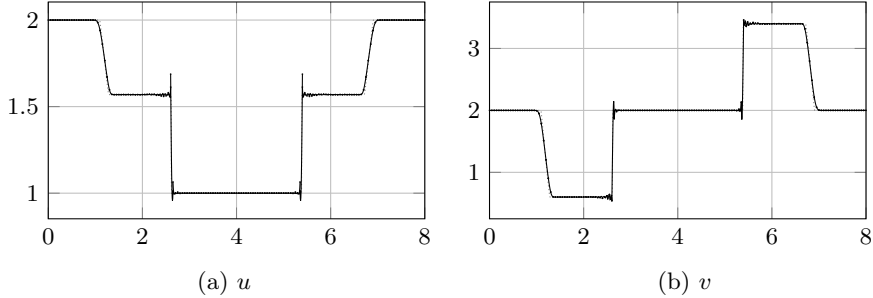

The amplitude of oscillations decreases as the degree of polynomials increases. For example, in Figure~\ref{fig:3} we consider cubic polynomials and $k/h=1/28$. 


It is worth mentioning that, compared to the standard discretization of the system as a conservation law using DG and Euler method in time, the optimization method has significantly less oscillations near discontinuities.  This does not hold though when a TVD Runge-Kutta (for example Osher's 3\textsuperscript{rd} order RK) time discretization is used. This can be seen in Table~\ref{tab:2}, where the exact value of total variation for $u$ is $2$ and for $v$ is $5.6$. 

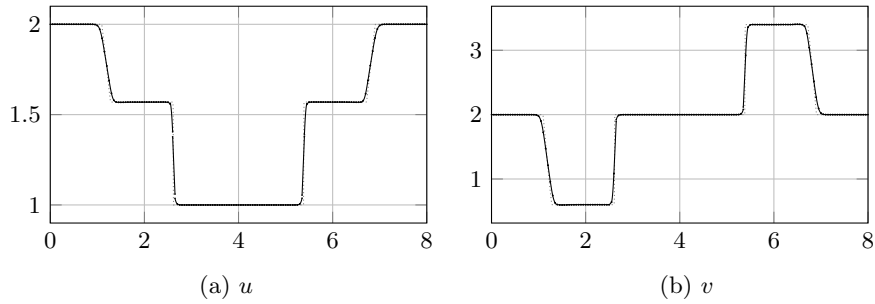
\begin{figure}[htbp] 
\centering
\begin{subfigure}{.49\textwidth} 
\begin{tikzpicture} 
\begin{axis}[discont axis]
\addplot[cont plot, gray, densely dotted] table[x=xn, y=u] {figs/discont_exact.txt};
\addplot[mark repeat=11, mark phase=10] table[x=x, y=u3_sl] {figs/discont_cubic.txt}; 
\addplot[only marks, mark repeat=11] table[x=x, y=u3_sl] {figs/discont_cubic.txt};
\end{axis}
\end{tikzpicture}
\caption{$u$}
\end{subfigure} 
\begin{subfigure}{.49\textwidth} 
\begin{tikzpicture} 
\begin{axis}[discont axis]
\addplot[cont plot, gray, densely dotted] table[x=xn, y=v] {figs/discont_exact.txt};
\addplot[mark repeat=11, mark phase=10] table[x=x, y=v3_sl] {figs/discont_cubic.txt}; 
\addplot[only marks, mark repeat=11] table[x=x, y=v3_sl] {figs/discont_cubic.txt};
\end{axis}
\end{tikzpicture}
\caption{$v$} 
\end{subfigure}
\caption{DG-IP, $N=160$, $T=0.25$, $k/h=1/28$, piecewise cubic polynomials, moments limiter. Discontinuous solution. \label{fig:4}} 
\end{figure}

\begin{table}[htbp] \ttfamily \small \centering 
\begin{tabular}{|c||c|c||c|c||c|c|} \hline 
~ & \multicolumn{2}{c||}{\rmfamily optimization} & \multicolumn{2}{c||}{\rmfamily DG Euler} & \multicolumn{2}{c|}{\rmfamily RKDG Osher} \\ \hline 
$N$ & $\operatorname{TV}u$ & $\operatorname{TV}v$ & $\operatorname{TV}u$ & $\operatorname{TV}v$ & $\operatorname{TV}u$ & $\operatorname{TV}v$ \\ \hline 
 40 & 2.269 & 6.601 & 2.627 & 7.301 & 2.369 & 6.509 \\ \hline 
 80 & 2.338 & 6.601 & 2.918 & 8.135 & 2.446 & 6.669 \\ \hline 
160 & 2.339 & 6.559 & 3.443 & 9.657 & 2.529 & 6.845 \\ \hline
320 & 2.294 & 6.416 & 3.858 & 11.00 & 2.522 & 6.883 \\ \hline
\end{tabular} 
\caption{DG-IP, LLF flux, comparison of total variation. $T=0.25$, $k/h=1/12$, $\lambda=0.25$, $\mu=1$. Discontinuous solution. Linear polynomials. \label{tab:2}}
\end{table}
To address the formation of oscillations we employ a slope limiting technique; 
depending on the degree of polynomials we select an appropriate method as described in Section~\ref{subsec:3p4}. 
The eigenvalues and left/right normalized eigenvectors required for the projection, for the current $\sigma$, are 
\[ 
\lambda_{1,2} = \pm \sqrt{3u^2+1},\ 
l_{1,2} = \frac 1 2\begin{pmatrix}  \frac{1}{\sqrt{3u^2 + 1}} & -\frac{1}{\sqrt{3u^2 + 1}}\\ 1 & 1\end{pmatrix},\ 
r_{1,2} = \begin{pmatrix} \phantom{+}\sqrt{3u^2 + 1}& 1\\ -\sqrt{3u^2 + 1} & 1\end{pmatrix}. 
\] 
The resulting solution for piecewise cubic polynomials can be seen in Figure~\ref{fig:4}. 

\FloatBarrier

\subsection{Optimization iteration convergence criteria}
\label{subsec:4p3} 

Suitable stop criteria for GD are necessary to avoid excessive iterations that will slow down the code and possibly introduce roundoff errors. 

\begin{figure}[htb] 
\centering
\begin{tikzpicture}[ 
every axis/.append style={x dir=reverse, xtick distance=1, width=.8\textwidth, height=.3\textwidth} 
]
\pgfplotsset{scale only axis}

\pgfplotstableread{
tol   err       iter
14 4.915e-03 78
13 4.915e-03 66
12 4.915e-03 58
11 4.915e-03 50
10 4.915e-03 42
09 4.907e-03 34
08 4.843e-03 26
07 4.922e-03 18
06 2.636e-02 10
05 1.963e-01 3
}\tdata

\begin{axis}[
  grid=major, 
  axis y line*=left, 
  xlabel=$-\log_{10} c_u$,
  ylabel=$\|u-u_\text{ex}\|$, 
  scaled y ticks=base 10:2,
  ytick distance=.05
]
\addplot[black,mark=x] table[x=tol, y=err] {\tdata}; 
\label{fig:5_error1}
\end{axis}
\begin{axis}[
  axis y line*=right,
  axis x line=none, 
  legend style={at={(.85,.95)},anchor=north east},
  ylabel=avg iterations
]
\addlegendimage{/pgfplots/refstyle=fig:5_error1, mark=x}
\addlegendentry{$\|u-u_\text{ex}\|$}
\addplot[densely dotted, mark=o, mark options=solid] table[x=tol, y=iter] {\tdata}; 
\addlegendentry{iterations}
\end{axis}
\end{tikzpicture}

\caption{GD iterations count and errors w.r.t.\ $c_u$. $T=0.25$, $N=160$, $k/h=1/12$, $\lambda=0.25$, $\mu=1$. Smooth solution. \label{fig:5} }
\end{figure}
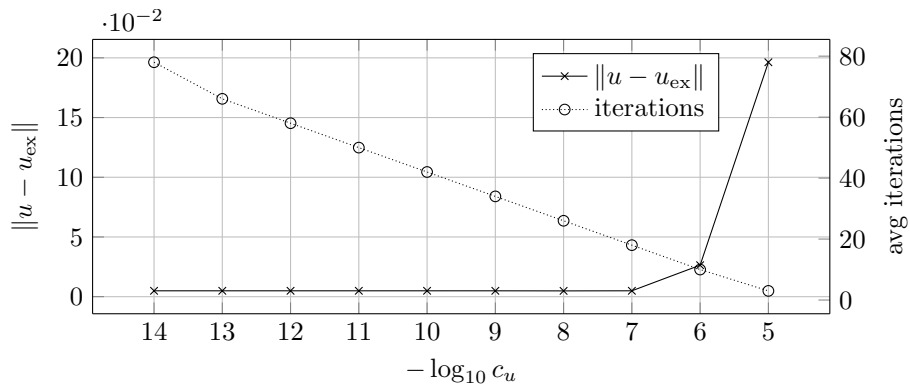

Consider for example the smooth initial conditions \eqref{eq:18}; errors and average GD iterations count per time step, in relation to convergence tolerances, can be seen 
in Figure~\ref{fig:5}. 
We see that $78$ iterations are required for the strictest tolerances, while only $14$ iterations are needed when we set $c_I=c_u=10^{-6}$ with insignificant increase in approximation error. 


\FloatBarrier 

\noindent Similar results hold for the discontinuous solution \eqref{eq:19}, as can be seen in Table~\ref{tab:3}. 

\begin{table}[htbp] \ttfamily \small \centering 
\begin{tabular}{|c|c||c||c|c|} \hline 
$\log_{10} c_I$ & $\log_{10} c_u$ & {\rmfamily avg iter} & $\|u-u_\text{ex}\|$ & $\|v-v_\text{ex}\|$    \\ \hline 
-14 & -14 & 86 & 1.192e-01 & 3.340e-01 \\ \hline
-10 & -10 & 55 & 1.192e-01 & 3.340e-01 \\ \hline
 -8 &  -8 & 39 & 1.192e-01 & 3.340e-01 \\ \hline
 -6 &  -6 & 23 & 1.191e-01 & 3.340e-01 \\ \hline
 -4 &  -4 &  9 & 1.177e-01 & 3.655e-01 \\ \hline
\end{tabular} 
\caption{GD iterations count for various stopping parameters. $T=0.25$, $N=80$, $k/h=1/12$, $\lambda=0.25$, $\mu=1$. Discontinuous solution.  \label{tab:3} }  
\end{table}

\FloatBarrier

\subsection{GD step selection} 
\label{subsec:4p4} 
GD step plays an important role in convergence speed of the method and can lead to significant acceleration if selected adaptively within each time step. Many methods exist in optimization literature for this purpose, such as e.g.  {\ttfamily adagrad}, {\ttfamily adam}, etc., see~\cite{KW19}) 
and may be adapted for our setting. 
Here, for illustrative purposes, we propose the following simple heuristic algorithm 
\begin{tabbing} 
\hspace{2em} \= \hspace{2em} \= \kill 
$\lambda_\text{min} = 1/4$, $\lambda_\text{lmax} = +\infty$, $c_\text{roff} = 1\mathrm e{-}10$ \\ 
$\lambda = \lambda_\text{init} = 1/4$ \\ 
for $l=1,\, 2,\, \ldots $ \+ \\ 
	calculate $u^j_{l+1}$, $v^j_{l+1}$ \\ 
	if $l \leq 2$ \\ 
		\> continue\\ 
	if $I_{l+1} < I_l$ \\ 
		\> $\lambda = \max\bigl(\frac 3 2 \lambda,\, \lambda_\text{lmax} \bigr)$\\ 
	else if $|I_{l+1}-I_l| > -5\,(I_l-I_{l-1})$ and $|I_{l+1}-I_l|>c_\text{roff}$\\ 
		\> $\lambda = \min\bigl(\frac 2 5 \lambda,\, \lambda_\text{min}\bigr)$\\ 
		\> $\lambda_\text{lmax} = \lambda$ \\ 
		\> $u^j_{l+1} = u^j_l$, $v^j_{l+1} = v^j_l$ 
\end{tabbing} 
where the $\lambda_\text{init}$ can be a reasonable initial value, or be chosen based on its value at previous time steps. 

\noindent As a typical example to highlight the performance benefits, 
we consider the smooth solution example, \eqref{eq:18}, with parameters $c_I=1\mathrm e{-}14$, $c_u=1\mathrm e{-}14$, $T=0.25$, $N=80$, $k/h=1/12$. This algorithm reduces the average number of required iterations from $73$ to $26$. 

Another significant measure for the performance of the algorithm is the approximation error given a fixed amount of computational resources. Using the aforementioned parameters we investigate the approximation error for various $N$ given that the maximum amount of iterations is limited, i.e.\ $c_i=10$. Errors, as well as the difference between the last two iterations of GD (that is $9$th and $10$th) of the integral under minimization and the value of $u$ can be seen in Table~\ref{tab:4}. 

\begin{table}[htbp] \ttfamily \small \centering 
\begin{tabular}{|c||c|c||c|c|} \hline 
$N$ & $\|u-u_\text{ex}\|$ & $\|v-v_\text{ex}\|$ & $|I_l-I_{l-1}|$ & $\|u_l - u_{l-1}\|$ \\ \hline \hline 
80 {\rmfamily ($\lambda$ \makebox[6.5ex][l]{fixed)}}  & 1.832e-02 & 1.072e-01 & 6.1e-07 & 5.5e-07 \\ \hline 
80 {\rmfamily ($\lambda$ \makebox[6.5ex][l]{adapt)}}  & 1.033e-02 & 2.664e-02 & 4.4e-10 & 1.8e-08 \\ \hline \hline 
320 {\rmfamily ($\lambda$ \makebox[6.5ex][l]{fixed)}} & 1.864e-02 & 9.989e-02 & 4.0e-08 & 4.0e-08 \\ \hline 
320 {\rmfamily ($\lambda$ \makebox[6.5ex][l]{adapt)}} & 2.355e-03 & 6.842e-03 & 1.5e-09 & 1.6e-08 \\ \hline 
\end{tabular} 
\caption{DG-IP, LLF flux, fixed versus adaptive $\lambda$. $T=0.25$, $c_i=10$. Smooth solution. Linear polynomials.  \label{tab:4} }
\end{table}

Finally we notice that the order of convergence is still maintained. We repeat the simulation of Section~\ref{subsec:4p1} for the smooth solution defined in \eqref{eq:18}. Convergence rates are as expected and can be seen in Table~\ref{tab:5}. 

\begin{table}[htbp] \ttfamily \small \centering 
\setlength{\tabcolsep}{4pt}
\begin{tabular}{|c||c|c|c|c||c|c|c|c|} \hline 
$N$ & $\|u-u_\text{ex}\|$ & {\rmfamily rate} & $\|u-u_\text{ex}\|_\infty$ & {\rmfamily rate} & $\|v-v_\text{ex}\|$ & {\rmfamily rate} & $\|v-v_\text{ex}\|_\infty$ & {\rmfamily rate} \\ \hline 
 20 & 3.164e-02 & -    & 4.141e-02 & -    & 1.061e-01 & -    & 1.132e-01 & -    \\ \hline 
 40 & 1.012e-02 & 1.65 & 1.638e-02 & 1.34 & 2.336e-02 & 2.18 & 3.049e-02 & 1.89 \\ \hline 
 80 & 2.564e-03 & 1.98 & 4.244e-03 & 1.95 & 5.332e-03 & 2.13 & 8.057e-03 & 1.92 \\ \hline 
160 & 6.279e-04 & 2.03 & 1.043e-03 & 2.02 & 1.300e-03 & 2.07 & 2.058e-03 & 1.97 \\ \hline 
320 & 1.554e-04 & 2.01 & 2.599e-04 & 2.01 & 3.225e-04 & 2.01 & 5.152e-04 & 2.00 \\ \hline 
\end{tabular} 
\caption{DG-IP, LLF flux, adaptive $\lambda$, convergence rates. $T=0.025$. Smooth solution. Linear polynomials.  \label{tab:5} }
\end{table}

\FloatBarrier ~


{\hfuzz=1.5em 
\ifxebiblatex 
\printbibliography[heading=bibintoc] 
\else 
\bibliography{ref}
\bibliographystyle{plain}
\fi
}

\end{document}